\documentclass[12pt,a4paper]{article}
\usepackage{amsmath,amssymb}
\newtheorem{thm}{Theorem}

\newtheorem{lem}{Lemma}
\newtheorem{prop}{Proposition}
\begin{document}
\title{Nef and big divisors \\
on   toric weak Fano 3-folds}
\author{Shoetsu OGATA\thanks{e-mail:  ogata{\char'100}math.tohoku.ac.jp}\\
Mathematical Institute, Tohoku University\\ Sendai 980-8578, Japan}
\date{October, 2013}
\maketitle
\begin{abstract}
We show that   a  nef and big line bundle whose adjoint bundle has non-zero
global sections on a nonsingular toric weak Fano 3-fold is normally generated.
As a consequence, we see that 
any ample line bundle on a nonsingular toric waek Fano 3-fold is normally generated.
As a corollary, we see that an ample line bundle 
whose adjoint bundle has non-zero global sections 
on a Gorenstein toric  Fano 3-fold is normally generated.
\end{abstract}

\section*{Introduction}

We call an invertible sheaf on an algebraic variety a line bundle.
A line bundle $L$ on an algebraic variety is called {\it normally generated}
(by Mumford\cite{Mf})
if the multiplication map of global sections $\Gamma(L)^{\otimes l}
\to \Gamma(L^{\otimes l})$ is surjective for all $l\ge1$.
We are interested in normal generation of ample line bundles on a toric variety.
If an ample line bundle $L$ on a normal algebraic variety $X$ is normally generated, then we 
see that it is very ample and that the graded ring $\bigoplus_{l\ge0} \Gamma(X, L^{\otimes l})$
is generated by elements of degree one and is a normal ring.
It is known that an ample line bundle on a nonsingular toric variety is always very ample
(see \cite[Corollary 2.15]{Ob}).  
We may ask whether any ample line bundle be normally generated.

In general, for an ample line bundle $L$ on a  (possibly singular) toric variety of dimension $n$,
we see that
\begin{equation}\label{int:1}
\Gamma(L^{\otimes l})\otimes \Gamma(L) \longrightarrow \Gamma(L^{\otimes(l+1)})
\end{equation}
is surjective for $l\ge n-1$ (see \cite{EW}, \cite{N} or \cite{NO}).  
When $n\le2$, hence, we see that all ample line bundles are 
normally generated (see \cite{K}).
We also have examples of  very ample but not normally generated line bundles for $n\ge3$
(see \cite{BG1}, \cite{BG2}, \cite{Og2} and \cite{Og3}).

An algebraic variety $X$ is called {\it Gorenstein} if its dualizing sheaf is invertible.
A Gorenstein variety $X$ is called {\it Fano} (or {\it weak Fano}) if
its anti-canonical divisor $-K_X$ is ample (or nef and big), respectively.
We know that the anti-canonical line bundle on a nonsingular
toric Fano variety of dimension $n$ is
normally generated if $n\le7$ (see \cite{HP}).
Ogata\cite{Og} shows that an ample line bundle $L$ on a nonsingular toric 3-fold $X$
 is normally generated if $h^0(L+K_X)=0$ or if $h^0(L+K_X)\not=0$ and $L+K_X$ is not big.

In this paper we restrict $X$ to be a nonsingular toric weak Fano 3-fold, but we assume 
that $L$ is a certain nef and big line bundle on $X$.

\begin{thm}\label{int:tm0}
Let $X$ be a nonsingular toric weak Fano variety of dimension three.
If a nef and big line bundle $L$ on $X$ satisties the condition that 
$h^0(L+K_X)\not=0$,
then $L$ is normally generated.
\end{thm}

By combining the result of \cite{Og} we obtain the following theorem.

\begin{thm}\label{int:tm01}
Let $X$ be a nonsingular toric weak Fano variety of dimension three.
Then any ample line bundle  on $X$  is normally generated.
\end{thm}

Ogata already shows that all ample line bundles on a nonsingular toric 3-fold with
non-trivail morphism onto the projective line are normally generated in \cite{Og4}.
Thus we have two classes of nonsingular toric 3-folds such that all
ample line bundles are normally generated.

For the proof of Theorem~\ref{int:tm0}, we need to prove the  weak version.

\begin{thm}\label{int:tm}
Let $X$ be a nonsingular toric weak Fano variety of dimension three.
If a nef and big line bundle $L$ on $X$ satisties the condition that $2L+K_X$ is nef, 
$h^0(L+K_X)\not=0$ and that $h^0(L+2K_X)=0$,
then $L$ is normally generated.
\end{thm}

We will give a proof of Theorem~\ref{int:tm} in Sections~\ref{sect:t1} and \ref{sect:t2}
by dividing into two cases.

Since a Gorenstein toric weak Fano 3-fold admits a crepant resolution, Theorem~\ref{int:tm0} 
implies the following theorem.

\begin{thm}\label{int:t2}
Let $Y$ be a Gorenstein toric weak Fano variety of dimension three.
If an ample line bundle $L$ on $Y$ satisfies the condition that $h^0(L+K_Y)\not=0$,
 then $L$ is normally generated.
\end{thm}

In our proof we do not use full classifications of Fano polytopes but use a classifycation of
minimal Fano polytopes of Kasprzyk\cite{Kp}.  There are 4,319 Gorenstein
toric Fano 3-folds (cf. \cite{KS}).

We note that there is an ample but not normally generated line bundle $L$ 
on a Gorenstein toric Fano 3-fold $Y$ 
with $h^0(L+K_Y)=0$.  See Remark in Section~\ref{sect6}.

\section{Line bundles on toric varieties}

In this section we recall the fact about toric varieties and line bundles on them from
Oda's book\cite{Ob} or Fulton's book\cite{Fb}.

Let $N$ be a free $\mathbb{Z}$-module of rank $n$ and $M:=\mbox{Hom}(N, \mathbb{Z})$
 its dual with the pairing
$\langle \cdot, \cdot\rangle : M\times N \to \mathbb{Z}$.
By scalar extension to $\mathbb{R}$, we have real vector spaces $N_{\mathbb{R}}:=
N\otimes_{\mathbb{Z}}\mathbb{R}$ and $M_{\mathbb{R}}:= M\otimes_{\mathbb{Z}}
\mathbb{R}$.  We also have the pairing of $M_{\mathbb{R}}$ and $N_{\mathbb{R}}$ by
scalar extension, which is denoted by the same symbol $\langle \cdot, \cdot\rangle$.
 
The group ring $\mathbb{C}[M]$ defines an algebraic torus $T_N:=\mbox{Spec}\ 
\mathbb{C}[M]\cong (\mathbb{C}^*)^n$ of dimension $n$.
Then the character group $\mbox{Hom}_{\mbox{gr}}(T_N, \mathbb{C}^*)$ of the algebraic
torus $T_N$ coincides with $M$.  For $m\in M$ we denote the corresponding
character by $e(m): T_N \to \mathbb{C}^*$.

Let $\Delta$ be a finite complete fan of $N$.  A convex cone $\sigma \in \Delta$ defines an affine
variety $U_{\sigma}:=\mbox{Spec}\ \mathbb{C}[M\cap \sigma^{\vee}]$.
Here $\sigma^{\vee}:= \{y\in M_{\mathbb{R}}; \langle y, x\rangle \ge0
\quad\mbox{for all}\quad x\in \sigma\}$ is the dual cone of $\sigma$.
Then we obtain a normal algebraic variety $X(\Delta) := \bigcup_{\sigma\in \Delta}
U_{\sigma}$, which is called a {\it toric} variety.
We note that $U_{\{0\}}\cong T_N$ is a unique dense $T_N$-orbit in $X(\Delta)$.
Set $\Delta(i):=\{\sigma\in \Delta; \dim \sigma =i\}$.
Then an element $\sigma\in \Delta(i)$ corresponds to a $T_N$-invariant subvariety
$V(\sigma)$ of dimension $n-i$.  In particular, $\Delta(1)$ corresponds to the set of 
all irreducible $T_N$-invariant divisors on $X(\Delta)$.

Let $\Delta(1)=\{\rho_1, \dots, \rho_s\}$ and $v_i$ the generator of the semi-group
$\rho_i\cap N$.  We simply write as $X=X(\Delta)$ and $D_i:= V(\rho_i)$ for $i=1, \dots, s$.
For a $T_N$-invariant line bundle $L$ there exists a $T_N$-invariant divisor
$D=\sum_i a_i D_i$ satisfying $L\cong \mathcal{O}_X(D)$.
For a $T_N$-invariant Cartier divisor $D$ we defines an rational convex polytope
$P_D \subset M_{\mathbb{R}}$ as
\begin{equation}\label{sect1:1}
P_D :=\{ y\in M_{\mathbb{R}}; \langle y, v_i\rangle \ge -a_i\quad \mbox{for}
\quad i=1, \dots, s\}.
\end{equation}
By definition we note that $P_{lD} =lP_D$ for any positive integer $l$.
Moreover, for another $T_N$-invariant Cartier divisor $E$ we have
$P_{D+E}=P_D +P_E$.  Here $P_D+P_E:=\{x+y\in M_{\mathbb{R}};
x\in P_D \quad\mbox{and}\quad y\in P_E\}$ is the Minkowski sum of $P_D$ and $P_E$.
By using this polytope, we can describe the space of global sections (see \cite[Section 2.2]{Ob},
or \cite[Section 3.5]{Fb})
\begin{equation}\label{sect1:2}
\Gamma(X, \mathcal{O}_X(D)) \cong \bigoplus_{m\in P_D\cap M} \mathbb{C}e(m).
\end{equation}
Hence, we see that  the surjectivity of the multiplication map of global sections
\begin{equation}\label{sect1:3}
\Gamma(X, \mathcal{O}_X(D)) \otimes \Gamma(X, \mathcal{O}_X(E))
\longrightarrow \Gamma(X, \mathcal{O}_X(D+E))
\end{equation}
is equivalent to  the equality
\begin{equation}\label{sect1:4}
P_D \cap M +P_E \cap M = (P_D+P_E) \cap M.
\end{equation}

If $\mathcal{O}_X(D)$ is generated by global sections, then all vertices of $P_D$ are
lattice points, that is, $P_D$ is the convex hull of finite subset of $M$.
Conversely, if for all $\sigma\in \Delta$ there exist $u(\sigma) \in M$ with
\begin{equation}\label{sect1:5}
\langle u(\sigma), v_i\rangle = -a_i \quad \mbox{for} \quad v_i\in \sigma
\end{equation}
and if $P_D$ is the convex hull of $\{u(\sigma); \sigma\in \Delta\}$,
then $\mathcal{O}_X(D)$ is generated by global sections (see \cite[Theorem 2.7]{Ob},
or \cite[Section 3.4]{Fb}).
We also knows \cite{Mv} that if $\mathcal{O}_X(D)$ is generated by global sections, then
there exists an equivariant surjective morphism $\pi: X \to Y$ to a toric variety $Y$ and an
ample line bundle $A$ on $Y$ with $\mathcal{O}_X(D)\cong \pi^*A$.
Thus we see that $\mathcal{O}_X(D)$ is generated by global sections if and only if
$D$ is {\it nef} (see also \cite[Theorem 3.1]{Ms}).

If $X$ is Gorenstein, then $-K_X=\sum_i D_i$ is a Cartier divisor.
By definition $P_{-K_X}$ is an integral polytope of dimension $n$ since the polytope is the 
intersection of half-spaces containing the origin as their interiors.
This implies that $-K_X$ is  {\it big}.

Now we introduce a criterion of nef-ness on nonsingular toric surfaces.

\begin{prop}\label{sect1:p1}
Let $X$ be a nonsingular complete toric surface and let $D$ a $T_N$-invariant divisor
with $|D|\not=\emptyset$.  If $|D|$ has no fixed components, then it is free from base points.
\end{prop}
{\it Proof}.
Since $\Delta(1)=\{\rho_1, \dots, \rho_s\}$ consists of half-lines from the origin in the plane
$N_{\mathbb{R}}$,
we may assume that $\rho_i$ and $\rho_{i+1}$ sit next to each other (as
usual we consider as $\rho_{s+1} =\rho_0$). 
Set $\sigma_i=\rho_i +\rho_{i+1}\in \Delta(2)$ for $i=1, \dots, s$.
Take $D = 
\sum_i a_i D_i$ with $|D|\not=\emptyset$. We may assume that $a_i\ge0$ for all $i$.

First we consider the case that $P_D$ is an integral convex polytope, that is, it is the 
convex hull of a finite subset of $M$.
Set $H^+(a_i):=\{y\in M_{\mathbb{R}}; \langle y, v_i\rangle \ge -a_i\}$ the half-plane
and its boundary line $H(a_i)$.
By definition~(\ref{sect1:1}) we see that $P_D$ is the intersection of all half-planes
$H^+(a_i)$'s.
Let $u_0$ be a vertex of $P_D$.  If $\dim P_D=2$, then a $1$-dimensional face of 
$P_D$ containing $u_0$ is contained in some line $H(a_i)$.  If $\dim P_D\le 1$,
then $P_D$ itself is contained in some $H(a_i)$.   We may set $i=1$.

Since $P_D$ is the intersection of $H^+(a_i)$'s, we take another line $H(a_j)$ ($j\not=1$)
meeting with $H(a_1)$ at $u_0$.
We may assume that all $\sigma_i$ with $i=1, \dots, j-1$ does not contain $-v_1$.
We claim that the line $H(a_i)$ contains $u_0$ for $i=2, \dots, j$.

For $\sigma_i =\rho_i + \rho_{i+1} \in \Delta(2)$,since $\{v_i,  v_{i+1}\}$ is a $\mathbb{Z}$-basis of $N$,
there exists $u(\sigma_i)\in M$ satisfying the condition~(\ref{sect1:5}).
Then we have
$$
u_0 \in H^+(a_1) \cap H^+(a_j) \subset u(\sigma_i) + \sigma_i^{\vee}
$$
for $i=1, \dots, j-1$.
If $u(\sigma_1)\not= u_0$, then the half-plane $H^+(a_2-1)$ would contain $P_D$.
This implies that $D_2$ is a fixed component of $|D|$.
Then we see that $u(\sigma_1)=u_0$. Considering $v_3, \dots, v_j$
 successively, we see that $u(\sigma_i)=u_0$
for $i=1, \dots, j-1$.

When $\dim P_D=2$, since we can take $H(a_j)$ so that it contains a 1-dimensional
face of $P_D$, we see that the opposite vertex on the edge $H(a_j)\cap P_D$
coincides with $u(\sigma_j)$.

When $\dim P_D\le1$, the vector $-v_1$ coincides with some $v_k$ ($j<k$).
By the same argument, we see that $u(\sigma_i)=u_0$ for $i=j, \dots, k-1$.
And we see that $u(\sigma_k)$ is also a vertex of $P_D$. 
  Hence, $\mathcal{O}_X(D)$ is
generated by global sections.

Next we assume only that $P_D$ is a rational convex polytope.
We can choose a positive integer $l$ so large that $lP_D$ is an integral polytope.
Since $lP_D= P_{lD}$, the line bundle $\mathcal{O}_X(lD)$ is generated by global sections,
hence it is nef.  Then $D$ is nef.  On a toric variety, if $D$ is nef, then $\mathcal{O}_X(D)$
is generated by global sections.  \hfill $\Box$

\medskip

{\bf Remark}.
If $\dim X\ge 3$, then the same statement of Proposition~\ref{sect1:p1} does not hold.
We can easily construct counterexamples, as Professor Payne points out.

\section{Adjoint line bundles}

Let $\omega_X$ be the dualizing sheaf on a toric variety $X$.  If a $T_N$-invariant
Cartier divisor $D$ is ample, then we have (see \cite[Proposition 2.24]{Ob})
$$
\Gamma(X, \mathcal{O}_X(D)\otimes \omega_X) \cong \bigoplus_{m\in (\mbox{\scriptsize{
Int}}(P_D))\cap M} \mathbb{C} e(m).
$$
If we take a resolution $\pi: \tilde X \to X$ of singularities by a subdivision of $\Delta$, then
$L=\pi^*\mathcal{O}_X(D)$ is nef and big, and we have
$$
\Gamma(\tilde X, L+K_{\tilde X}) \cong \Gamma(X, \mathcal{O}_X(D)\otimes\omega_X).
$$

In \cite{Og} we show that an ample line bundle $A$ on a nonsingular toric 3-fold $X$
 is normally generated if $A$ satisfies the condition
  that $h^0(X, A+K_X)\not=0$ and that $A+K_X$ is not big.
In order to treat more general case, we have to know the adjoint bundle $L+K_X$
with $h^0(L+K_X)\not=0$ for a nef and big line bundle $L$.

\begin{lem}\label{sect2:l1}
Let $X$ be a nonsingular complete toric variety of dimension three.
If a nef and big line bundle $L$ on $X$ satisfies that $h^0(X, L+K_X)\not=0$, 
then the fixed part of $L+K_X$ is a reduced divisor $\sum_i E_i$ with
$(E_i, L_{E_i}) \cong (\mathbb{P}^2, \mathcal{O}_{\mathbb{P}^2}(1))$ and $\mathcal{O}_X(E_i)|_{E_i}\cong 
\mathcal{O}_{\mathbb{P}^2}(-1)$, and $L+K_X-\sum_i E_i$ is nef.
\end{lem}
{\it Proof}.
Since a nef and big line bundle on a toric variety is the pull back of an ample line bundle on
a toric variety of the same dimension, $K_X+tL$ is generated by global sections,
that is, nef  for a sufficiently large $t>0$.
If $K_X+tL$ is not nef, then Mori's theory \cite{Mr} says that
there exists an extremal curve $R$ with $(K_X+tL)R<0$ and $-K_XR\le 4$.
Thus we have $1\le t\le 3$.

By the Mori-Kawamata theory(cf. \cite{KMM}, \cite{Mr})
if $K_X+tL$ is not nef, then
we have a contraction  morphism $\varphi: X \to Y$. 
Following the same argument of Fujita \cite[Theorem 11.8]{Fj}, our assumption $h^0(X, L+K_X)\not=0$
implies that $t=1$ and $\varphi$ is birational and moreover
that $\varphi$ is a  blowing-ups of a point and 
there exists a nef and big line bundle $\bar L$
on $Y$ such that $L+E\cong \varphi^*\Bar L$, where $E$ is the
 exceptional divisor $E\cong \mathbb{P}^2$, $L_{E}\cong\mathcal{O}(1)$ and $\mathcal{O}_X(E)|_E\cong \mathcal{O}(-1)$.
Since $K_X=\varphi^*K_Y+2E$,  we have $L+K_X = \varphi^*(\bar L+K_Y) +E$.
\hfill  $\Box$
\medskip

\section{Proof of Theorem~\ref{int:tm} (Part I)}\label{sect:t1}

Let $B:=\sum_i D_i$ be the boundary divisor of $T_N$ in $X$.  We assume that 
$B$ is nef.  Let $L=\mathcal{O}_X(D)$ be a nef and big line bundle
satisfying the condition in Theorem~\ref{int:tm}.  Then
$P_D$ is an integral polytope of dimension three.
The assumption $h^0(X, L+K_X)\not=0$ of Theorem~\ref{int:tm} implies that
$(\mbox{Int}(P_D))\cap M \not= \emptyset$.
Let $F$ be the fixed components of $|D+K_X|$ and $A:=(D+K_X)-F$.  
From Lemma~\ref{sect2:l1} we see that $|A|$ is
free from base points.  
Since $\Gamma(X, L+K_X)=\Gamma(X, \mathcal{O}_X(A))$,
we see that $P_A$ coincides with the convex full of
$(\mbox{Int}(P_D))\cap M$.
We note that if $-K_X=B$ is nef, then $D-F=A+B$ is also nef.

Before treating nef divisors on 3-folds, we need to know more about
nef divisors on toric surfaces. For this purpose we heavily use the following lemma
given by Fakhruddin\cite{F}.

\begin{lem}\label{sect3:l0}
Let $I=[a, b] $ and $J=[c, d] $ be closed intervals in $\mathbb{R}$ with
$a, b\in \mathbb{Z}$ and $J\cap \mathbb{Z} \not=\emptyset$.
Then we have
$$
I\cap \mathbb{Z} + J\cap \mathbb{Z} = (I+J)\cap \mathbb{Z}.
$$
\end{lem}

By using this lemma we can prove the following lemma, which is an answer to the
Oda's question\cite{Od}.
\begin{lem}\label{sect3:l1}
Let $A$ and $B$ be nef divisors on a nonsingular complete toric surface $Y$.  Then
the multiplication map of global sections
$$
\Gamma(Y, \mathcal{O}_Y(A))\otimes \Gamma(Y, \mathcal{O}_Y(A+B))
\longrightarrow \Gamma(Y, \mathcal{O}_Y(2A+B))
$$
is surjective.
\end{lem}
{\it Proof}.
Since $\dim Y=2$, in this proof we set $M\cong \mathbb{Z}^2$ and $P_A, P_B \subset
M_{\mathbb{R}}\cong \mathbb{R}^2$.
We will show the equality
\begin{equation}\label{sect3:1}
P_A\cap M +(P_A+P_B) \cap M = (2P_A+P_B) \cap M.
\end{equation}

When $\dim P_A=1$, a generator of 
the sub-lattice $(\mathbb{R}P_A)\cap M\cong \mathbb{Z}$ is a part of a basis of $M$.
Let $\{u_1, u_2\}$ be a basis of $M$ such that $u_1$ is a generator of 
$(\mathbb{R}P_A)\cap M\cong \mathbb{Z}$.
By taking an affine transformation of $M$,
we may set $P_A=\{r u_1; 0\le r\le b\}$ for some positive integer $b$.
Then $P_A+P_B$ has two edges parallel to $P_A$.
Let $l_k:=\{xu_1+ku_2; x\in \mathbb{R}\}$ for an integer $k$.
Then we have a decomposition of lattice points in $2P_A+P_B$ as
$$
(2P_A+P_B)\cap M =\bigcup_{k\in \mathbb{Z}} (2P_A+P_B)\cap l_k\cap M.
$$
We can apply Lemma~\ref{sect3:l0} to line segments $(2P_A+P_B)\cap l_k$.

When $\dim P_A=2$, we can decompose $P_A$ into a union of basic triangles:
$$
P_A = \bigcup_i Q_i,
$$
where a basic triangle means that $Q_i\cap M$ has only three elements, that is,
vectors of two edges  are a generator of $M$.
Since $2P_A+P_B= \cup_i (2Q_i+P_B)$, we can reduce to the case that $P_A$ is
a basic triangle.  Moreover, when $\dim P_B=2$, we can also reduce to the case
that $P_B$ is a lattice triangle.   Since a parallel transformation of polytopes
by an element of $M$ does not change the equality(\ref{sect3:1}), we may assume that 
one vertex of polytopes is the origin. 
Set $E_1, E_2$ and $E_3$ be edges of $P_A$, and let $m_1, m_2$ and $0$ be
three vertices of $P_B$.  Then $2P_A+P_B$ is the union of $0+P_B$ and
$2E_i+P_B$ ($i=1,2,3$) and $2P_A+m_j$ ($j=1,2$)
(when $\dim P_B=1$ it is considered as $m_1=m_2$).
Then we can apply Lemma~\ref{sect3:l0}.
\hfill  $\Box$

\medskip

We return to the case that $B=-K_X$ and $A=D+K_X -F$ in dimension three.
First we prove Theorem~\ref{int:tm} in a special case.

\begin{prop}\label{sect3:p1}
When $\dim P_A\le2$ the bundle $\mathcal{O}_X(A+B)$ is normally generated.
\end{prop}
{\it Proof}.
Set $L=\mathcal{O}_X(A+B)$.  Since $\mathcal{O}_X(A)=L+K_X$, we have an
exact sequence
\begin{equation}\label{sect3:2}
0 \to \mathcal{O}_X(A) \to L \to L_B \to 0.
\end{equation}
Since $A$ is nef, we have $H^i(X, \mathcal{O}_X(A))=0$ for $i\ge1$.  Thus the
global sections of (\ref{sect3:2}) is exact.  Take the tensor product with
$\Gamma(X, \mathcal{O}_X(A))$.
When $\dim P_A\le2$, we see that $\mathcal{O}_X(A)$ is normally generated (see
(\ref{int:1})).

On the other hand, $\Gamma(B, (2L+K_X)_B)$ has a basis $\{e(m); m\in (\partial
(2P_A+P_B))\cap M\}$ as vector spaces.
One $e(m)$ is contained in $\Gamma(D_i, (2L+K_X)_{D_i})$ for some $D_i$.
Since the restriction map $\Gamma(X, G) \to \Gamma(D_i, G_{D_i})$ is surjective
for any nef line bundle $G$ on a toric variety $X$, from Lemma~\ref{sect3:l1} we see
that the multiplication map
$$
\Gamma(B, L_B)\otimes \Gamma(B, (L+K_X)_B) \longrightarrow
\Gamma(B, (2L+K_X)_B)
$$
is surjective.  Thus we obtain the surjectivity of
$\Gamma(L)\otimes \Gamma(L+K_X) \to \Gamma(2L+K_X)$.

By tracing the same argument after changing $A$ with $L=A+B$, we obtain
the normal generation of $\mathcal{O}_X(A+B)$.
\hfill  $\Box$

\section{Fano polytopes}

For our proof of Theorem, we need
 some properties of Fano polytopes of dimension three.
Batyrev \cite{B} called {\it Fano polytope} a lattice polytope containing 
no lattice points except the origin
in its interior.
Kasprzyk \cite{Kp} called a Fano polytope $P$ of dimension $n$ {\it minimal} if
 the convex hull $\mbox{Conv}\{P\cap (\mathbb{Z}^n)\setminus \{v\}\}$
is never Fano for any vertex $v$ of $P$.  He classifies minimal Fano polytopes
of dimension three.

Let $Y$ be a Gorenstein toric  Fano 3-fold and $R\subset M_{\mathbb{R}}$ the
lattice polytope corresponding to the ample anti-canonical divisor $-K_Y$.
$R$ contains only one lattice points in its interior.  By an affine transformation of $M$,
we may assume the origin is the interior point of $R$.  Then $R$ is a Fano polytope.
Let $F$ be a facet of $R$ and $G$ the minimal lattice triangle  with vertices in
$F\cap M$.  By a suitable choice of the coordinates of $M_{\mathbb{R}}$,
we may set $G=\{(1,0,0),(0,1,0),(0,0,1)\}$.
We fix the coordinates of $M_{\mathbb{R}}$ in this section.

\begin{lem}[Kasprzyk]\label{4:l1}
Let $S:=\mbox{\rm Conv}\{G, (-a,-b,-c)\}$ be a minimal Fano polytope with $1\le a\le b\le c$.
Then the triplet $(a,b,c)$ coinsides with one of
\begin{eqnarray*}
& &(1,1,1), (1,1,2), (1,1,3), (1,2,2), (1,2,3), (1,2,4)\\
& &(1,3,4), (1,3,5), (1,4,6), (2,3,5), (3,4,5).
\end{eqnarray*}
\end{lem}

We note that $S$ contains all $(\frac12,0,0),(0,\frac12,0),(0,0,\frac12), (\frac12,\frac12,0),(\frac12,0,\frac12), (0,\frac12,\frac12)$ in $\frac12 M$.
Here we call the point $v\in \frac12 M$ {\it perfect half-integral} if its any coordinates
are not integers.

\begin{lem}\label{4:l2}
$S$ contains a perfect half-integral point unless $(a,b,c)=(1,2,4)$.
\end{lem}

Set $S_0:=\mbox{Conv}\{G, (-1,-2,-4)\}$.  We note $(0,0,-1)\in S_0$.

\begin{lem}\label{4:l3}
A minimal Fano polytope of dimension three such that it contains $G$ as a facet
and has vertices more than four is one of 
\begin{itemize}
\item[{\rm(1)}] $S_1:=\mbox{\rm Conv}\{G, (0,0,-1),(-1,-1,0)\}$,
\item[{\rm(2)}] $S_2:=\mbox{\rm Conv}\{G, (0,0,-1),(-2,-1,0)\}$,
\item[{\rm(3)}] $S_3:=\mbox{\rm Conv}\{G, (1,1,-1),(-1,-1,0)\}$,
\item[{\rm(4)}] $S_4:=\mbox{\rm Conv}\{G, (-1,0,-1),(-2,-1,0)\}$,
\item[{\rm(5)}] $S_5:=\mbox{\rm Conv}\{G, (-1,0,-2),(0,-1,-2)\}$,
\item[{\rm(6)}] $S_6:=\mbox{\rm Conv}\{G, (0,0,-1),(0,-1,0),(-1,0,0)\}$.
\end{itemize}
\end{lem}

\begin{lem}\label{4:l4}
$S_1, S_3, S_4$ contain a perfect half-integral point.
\end{lem}

We note that $S_2, S_5, S_6$ contain the point $(0,0,-1)$.

\begin{prop}\label{4:p1}
Let $Y$ be a Gorenstein toric Fano 3-fold and $R$ the lattice polytope
corresponding $-K_Y$.  If $Q \subset M_{\mathbb{R}}$ is a lattice polytope containing
$\frac{m}2+\frac32 R$, then $\mbox{\rm Int}(Q)\cap M\not=\emptyset$.
\end{prop}
{\it Proof}.
 If
\begin{equation}\label{4:e1}
\left\{\frac{m}2 +\mbox{\rm Int}\left(\frac32 R\right)\right\}\cap M \not=\emptyset,
\end{equation}
then $Q$ contains lattice points in its interior.

If $\frac{m}2\in M$, then $R\subset \mbox{Int}(\frac32R)$ containes the origin.

Set $\frac{m}2 \notin M$. 
Then (\ref{4:e1}) is equivalent to the existence an element
\begin{equation}\label{4:e2} 
\frac{u}2 \in
\left\{\mbox{\rm Int}\left(\frac32 R\right)\right\}\cap \left\{\frac12M\setminus M\right\} 
\quad\mbox{with}\ \frac{m}2+\frac{u}2\in M.
\end{equation}

We may replace $R$ with a minimal Fano polytope $S$.
We fix the coordinates of $M_{\mathbb{R}}$ as $G=\{(1,0,0),(0,1,0),(0,0,1)\}$
for a minimal lattice triangle $G$ contained in a facet of $R$.
From Lemmas~\ref{4:l1}, \ref{4:l2}, \ref{4:l3} and \ref{4:l4},
we see that (\ref{4:e2}) holds for a minimal Fano polytope $S$ 
unless $\frac{m}2$ is a perfect half-integral point and $S$ is one of $S_0, S_2, S_5$ and $S_6$.

We assume that $\frac{m}2$ is a perfect half-integral point and $S=S_0, S_2, S_5$, or $S_6$.
Since $S$ contains the point $(0,0,-1)$, 
two half-integer points $(\frac12,\frac12,\pm\frac12)$ are contained in the boundary of
$\frac32 S$.
Set $\frac{m}2=(\frac12,\frac12,\frac12)$ and $S_m:=\frac{m}2 +\frac32 S$ the rational convex
polytope.   Then $S_m$ contains the lattice points $(1,1,1),(1,1,0)$.
Let $F_1, F_2$ be the facets of $S_m$ containing $(1,1,1)$ and $(1,1,0)$, respectively.
Then the intersection $F_1\cap F_2$ is the line segment $\ell$ conecting $(2, \frac12,\frac12)$ and
$(\frac12, 2,\frac12)$.
Let $Q$ be a lattice polytope containing $\frac{m}2+\frac32S$.
Assume that $Q$ contains no lattice points in its interior.
Since $(1,1,1),(1,1,0)$ are contained in the boundary of $Q$ and since these lattice points are
contained in the interior of the facets $F_1$ or $F_2$,  $Q$ has two facets $G_1, G_2$ containing
$F_1, F_2$, respectively.  In particular, the intersection $G_1\cap G_2$  contains $\ell$.
Thus $Q$ has a face $E$ of dimension one containing $\ell$.   Since $Q$ is a lattice polytope,
the end points of $E$ are lattice points.  But it is impossible.
Hence, we have $(\mbox{Int}\ Q)\cap M\not=\emptyset$.  \hfil $\Box$

\section{Proof of Theorem~\ref{int:tm} (Part II)}\label{sect:t2}

We recall our situation.  Let $X$ be a nonsingular toric weak Fano 3-fold and $B=\sum_i D_i$
the nef anti-canonical divisor.  Let $A$ be a nef and big divisor with 
$h^0(X, \mathcal{O}_X(A)+K_X)
=0$.  We consider the normal generation of $ \mathcal{O}_X(A+B))$.

\begin{prop}\label{sect3:p2}
If $\dim P_A=3$ and if $(\mbox{\rm{Int}}(P_A))\cap M=\emptyset$,
then $\mathcal{O}_X(A+B)$ is normally generated.
\end{prop}
{\it Proof}.
Assumption $(\mbox{\rm{Int}}(P_A))\cap M=\emptyset$ implies that 
$h^0(X, \mathcal{O}_X(A+K_X))=0$.

First we consider the multiplication with $\Gamma(\mathcal{O}_X(B))$.
Since $\mathcal{O}_X(B+K_X)\cong \mathcal{O}_X$, we have an exact sequence
$$
0\to \mathcal{O}_X \to \mathcal{O}_X(B) \to \mathcal{O}_B(B|_B) \to 0.
$$
Since the multiplication map $\Gamma(D_i, \mathcal{O}_{D_i}((A+B)|_{D_i}))
\otimes \Gamma(D_i, \mathcal{O}_{D_i}(B|_{D_i})) \to \Gamma(D_i, \mathcal{O}_{D_i}((A+2B)|_{D_i}))$ is surjective from Lemma~\ref{sect3:l1}, we have the surjectivity of
$\Gamma(X, \mathcal{O}_{X}(A+B))
\otimes \Gamma(X, \mathcal{O}_{X}(B)) \to \Gamma(X, \mathcal{O}_{X}(A+2B))$.

Set $L=\mathcal{O}_X(A+2B)$.  
Next we consider the multiplication of $\Gamma(L)$ and $\Gamma(\mathcal{O}_X(A))$.
In order that, we will show the vanishing of
$H^i(X, L(-iA))$ for $i\ge1$.

We have $H^1(X, L(-A)) = H^1(X, \mathcal{O}_X(2B))=0$ since $B$ is nef.

From the Serre duality  we have
$h^3(X, L(-3A))=h^3(X, \mathcal{O}_X(2B-2A))=h^0(X, \mathcal{O}_X(2A+3K_X))$.
If $h^0(X, \mathcal{O}_X(2A+3K_X))\not=0$, we have an injective homomorphism
$\mathcal{O}_X(3B) \to \mathcal{O}_X(2A)$.
By taking global sections, we have an inclusion
$$
m+P_{3B} \subset P_{2A}
$$
for some $m\in M$.  We may write it as
\begin{equation}\label{5:e1}
\frac{m}2 +\frac32 P_B \subset P_A.
\end{equation}
This implies that $P_A$ has lattice points in its interior from Proposition~\ref{4:p1}.
This contradicts to the assumption $(\mbox{\rm{Int}}(P_A))\cap M=\emptyset$.
Hence we have $h^3(L(-3A))=0$.

Consider the exact sequence
\begin{equation}\label{sect3:3}
0 \to \mathcal{O}_X(A+2K_X) \to \mathcal{O}_X(A+K_X) \to
\mathcal{O}_B((A+K_X)_B) \to 0.
\end{equation}
Since $H^j(X, \mathcal{O}_X(A+K_X))=0$ for $j=0,1$, if $h^1(X, \mathcal{O}_X(A+2K_X))\not=0$,
then $h^0(B, \mathcal{O}_B( ( A+K_X)_B)\not=0$.
A non-zero section defines an injective homomorphism
$\Gamma(D_i, \mathcal{O}_{D_i}(B|_{D_i})) \to \Gamma(D_i, \mathcal{O}_{D_i}(A|_{D_i}))$
for each irreducible invariant divisor $D_i$.
Since all facets of $P_B$ are contained in $P_A$, $P_B\subset P_A$ and $(\mbox{Int}P_A)
\cap M\not=\emptyset$.  This contradicts the assumption on $A$.
thus we have $h^1(X, \mathcal{O}_X(A+2K_X))=0$.

Considering the exact sequence
$$
0 \to \mathcal{O}_X(A+3K_X) \to \mathcal{O}_X(A+2K_X) \to
\mathcal{O}_B((A+2K_X)_B) \to 0,
$$
we can show $h^1(X, \mathcal{O}_X(A+3K_X))=0$ in the same way.
By Serre duality, we have $h^2(X, \mathcal{O}_X(-2K_X-A))=h^2(X, L(-2A))=0$.

From  vanishing of $H^i(X, L(-iA))$ for $i\ge1$ we can apply \cite[Theorem 2]{Mf}
to obtain the surjectivity of the multiplication map
\begin{equation}\label{sect3:4}
\Gamma(X, \mathcal{O}_X(A)) \otimes \Gamma(X, L) \longrightarrow
\Gamma(X, L(A)).
\end{equation}
Since $L=\mathcal{O}_X(A+2B)$, 
 we obtain the normal generation of $\mathcal{O}_X(A+B)$
from the first step of the proof and the surjectivity of (\ref{sect3:4}).
\hfill  $\Box$

\medskip

Finally we will complete the proof of Theorem~\ref{int:tm}.

{\it Proof of Theorem~\ref{int:tm}}.
By assumption that $L$ is  nef and big with $h^0(X, L+K_X)\not=0$
and $h^0(L+2K_X)=0$.

If $L+K_X$ has no fixed components, then we see the normal generation of $L$
from Propositions~\ref{sect3:p1} and \ref{sect3:p2}.
Let $F$ be the fixed components of $L+K_X$.  Then $F=\sum_i E_i$,$E_i\cong
\mathbb{P}^2$ and $E_i$'s are disjoint from Lemma~\ref{sect2:l1}.
And we have $L_{E_i}\cong \mathcal{O}_{\mathbb{P}^2}(1)$ and
$L(-F)_{E_i}\cong \mathcal{O}_{\mathbb{P}^2}(2)$.

Consider the exact sequence
\begin{equation}\label{sect3:5}
0 \to L(-F) \to L \to L_F \to 0.
\end{equation}
Since $L(-F)$ is nef, we have $H^1(X, L(-F))=0$.  Thus the sequence of global 
sections of (\ref{sect3:5}) is exact.  Taking the tensor product with $\Gamma(X, L(-F))$,
we see the surjectivity of the map
$$
\Gamma(X, L(-F))\otimes \Gamma(X, L) \longrightarrow \Gamma(X, 2L(-F))
$$
since $L(-F)$ is normally generated  from Propositions~\ref{sect3:p1}
and \ref{sect3:p2}.
By changing the role of $\Gamma(X, L(-F))$ with $\Gamma(X,L)$ 
we see the normal generation of $L$.
\hfill $\Box$

\section{Proof of Theorems}\label{sect6}

{\it Proof of Theorem~\ref{int:tm0}}. 
Let $L$ be a nef and big line bundle on a nonsingular toric weak Fano 3-fold $X$ with
$h^0(X, L+K_X)\not=0$.  Set $F_1$ the fixed part of $L+K_X$.  Then $L+K_X-F_1$ is nef
from Lemma~\ref{sect2:l1}.    Set $L_1:=L+K_X-F_1$.
Set $P_1:=P_{L_1}$ the integral convex polytope corresponding to the nef line bundle $L_1$.
If $\dim P_1\le2$ or if $\dim P_1=3$ and $\mbox{Int}(P_1)=\emptyset$, then
$L$ is normally generated from Theorem~\ref{int:tm}.

Assume that $\dim P_1=3$ and $\mbox{Int}(P_1)\not=\emptyset$.
Then $L_1$ is nef and big and $h^0(X, L_1+K_X)\not=0$.
Set $F_2$ the fixed part of $L_1+K_X$.  Then $L_2:=L_1+K_X-F_2$ is nef.
If $h^0(X, L_2+K_X)=0$, then $L_1$ is normally generated from Theorem~\ref{int:tm}.
Set $B=\sum_iD_i$ the boundary divisor.  Since $L-F_1$ is a sum of two nef and big 
line bundles $L_1$ and $-K_X$, the fact that the short exact sequence
$$
0\to L_1\to L_1-K_X \to (L_1-K_X)|_B \to 0
$$
gives the short exact sequence of its global sections and 
the surjectivity of two multiplication maps
$\Gamma(L_1)\otimes\Gamma(L_1)\to \Gamma(2L_1)$ and
$\Gamma(B, (L_1-K_X)|_B)\otimes\Gamma(B, L_1|_B)\to \Gamma(B, (2L_1-K_X)|_B)$
implies the surjectivity of the spaces of global sections of $L_1-K_X$ and $L_1$.
By the same argument of the proof of Proposition~\ref{sect3:p2}, we have the surjectivity
of $2L_1-K_X$ and $-K_X$.  This implies the normal generation of $L-F_1=L_1-K_X$.
Next we apply the same argument of "Proof of Theorem~\ref{int:tm}"  in the previous section
to obtain the normal generation of $L$. 

If $h^0(X, L_2+K_X)\not=0$, then we continue the same argument for $L_2$.  By induction we 
obtain the proof of Theorem~\ref{int:tm0}.
\hfill $\Box$

\bigskip

If the anti-canonical divisor $-K_X$ of a Gorenstein toric variety $X$ is nef, then
it is nef and big, hence, there exists a polarized toric variety $(Y, A)$ and
a surjective morphism $\pi: X \to Y$ such that $-K_X\cong \pi^*A$.
Since $Y$ has only rational singularity, we see that $A= -K_Y$ and $Y$ is a 
Gorenstein Fano 3-fold.

On the other hand, let $Y$ be a Gorenstein toric weak Fano 3-fold.  Then we have a resolution
of singularities $\pi: X \to Y$ with $K_X\cong \pi^*K_Y$ (a crepant resolution).  
Thus we can apply Theorem~\ref{int:tm}
to a nef and big line bundle $\pi^*A$ with ample $A$ on $Y$.
We obtain Theorem~\ref{int:t2}.

\medskip

{\bf Remark}.
In Theorem~\ref{int:tm} or \ref{int:t2} we cannot remove the condition $h^0(X, L+K_X)\not=0$.
We give an example of $(X, L)$ such that $-K_X$ is nef but $L$ is not normally
generated and $h^0(X, L+K_X)=0$.

Let $M=\mathbb{Z}^3$ and $P:=\mbox{Conv}\{0, (1,0,0), (0,1,0), (1,1,2)\}$
in $M_{\mathbb{R}}$.
Then there exists the polarized toric 3-fold $(Y, \mathcal{O}_Y(D))$ with
$P_D=P$.  This $Y$ is Gorenstein toric Fano with $-K_Y =2D$.
Since $P$ does not contain lattice points of the form $(a,b,1)$,
we can easily see that  $D$ is not very ample.
We can make a toric crepant (partial) resolution $\pi: X \to Y$ of singularities with $K_X=\pi^*K_Y$.
Then $-K_X$ is nef (and big) and $L:=\pi^*\mathcal{O}_Y(D)$ is nef and big, and
$h^0(X, L+K_X)=0$.
We note that a chice of (partial) resolution $\pi: X \to Y$ is not unique for this Gorenstein toric
Fano 3-fold $Y$.

\end{document}